\documentclass[12pt]{article}
   
\newcommand{\p}{\partial}
\newtheorem{theorem}{Theorem}

\begin{document}

\title{A Remark on the Three Dimensional Baroclinic 
        	Quasigeostrophic Dynamics 
\footnote{This research was supported by the National Science
Foundation Grant DMS-9704345  and a Clemson
University Research Grant.
This research was performed while the 
author was visiting the Isaac Newton Institute for 
Mathematical Sciences at Cambridge University, England.}  }

\author{Jinqiao Duan \\
Department of Mathematical Sciences,\\
Clemson University, Clemson, South Carolina 29634, USA.\\
E-mail: duan@math.clemson.edu, Fax: (864)656-5230.} 

\date{January 12, 1998  }

\maketitle

\begin{abstract}

We consider time-periodic patterns of the dissipative
three dimensional baroclinic quasigeostrophic model 
in spherical coordinates, under time-dependent forcing. 
We show that when the forcing is time-periodic and
the spatial square-integral of the  forcing is bounded
in time, the model has time-periodic solutions.
 
\bigskip
{\bf Keywords---} quasigeostrophic fluid model,  dissipative dynamics, 
	          time-periodic motion,  nonlinear analysis.

\end{abstract}

\bigskip
\bigskip

{\bf Short running title:}

	Three Dimensional Quasigeostrophic Dynamics

\newpage

\section{Introduction}

The three dimensional baroclinic quasigeostrophic  model  
is an approximation of  the rotating
Euler (inviscid case) or Navier-Stokes (viscous case)
equations in the limit of zero Rossby number,
i.e., at asymptotically high rotation rate; see, for example, 
\cite{Charney}, \cite{Pedlosky}, \cite{Gill}, \cite{Beale},  
\cite{Embid_Majda}, \cite{Desjardins1}, and \cite{Desjardins2}.    

The well-posedness for the invisid three dimensional 
baroclinic quasigeostrophic  model
was studied in, for example, \cite{Dutton} and \cite{Beale}, and
for the viscous three dimensional baroclinic quasigeostrophic  
model in \cite{Bennett}.
The invisid baroclinic quasigeostrophic  model is a 
Hamiltonian system (\cite{Holm}).
The work mentioned above and most other
research about the baroclinic quasigeostrophic model has been  on
the Cartesian, $\beta-$plane version of this model.

Wang (\cite{Wang}) discussed global attractors for
the viscous three dimensional baroclinic quasigeostrophic  
model in sphericcal coordinates,
under time-independent forcing.

In this paper, we  consider the same
viscous three dimensional 
baroclinic quasigeostrophic  model, but under time-dependent forcing.
We show  that when the forcing is time-periodic and 
when the spatial square-integral of the forcing is bounded
in time, the forced viscous three dimensional
quasigeostrophic  model has time-periodic solutions. 
We use a topological technique from
nonlinear global analysis (\cite{Krasnoselskii}).

\section{Dissipative Dynamics and Time-Periodic Motion}

We consider the three dimensional baroclinic quasigeostrophic 
equation, in nondimensional form,  for the atmospheric dynamics 
in spherical coordinates  $\phi, \theta, \zeta$ as in \cite{Wang} 
\begin{eqnarray}
(A \psi)_t +J(\psi, Ro\; A\psi + 2cos\theta) - \frac1{Re}A^2\psi
	& = & f(\phi, \theta, \zeta, t) , \label{eqn}  
\end{eqnarray} 
where $\psi(\phi, \theta, \zeta, t)$ is the stream function,
$f(\phi, \theta, \zeta, t)$ is the time-dependent
forcing such as external source of heating,
and $Ro, Re$ are the Rossby number, Reynolds number, respectively.  
The space domain for the equation is $D=S^2 \times (\zeta_0, 1)$
with $S^2$ the two dimensional unit sphere and
$\zeta_0$ fixed: $1> \zeta_0 \geq 0$.
The operators appeared in the equation (\ref{eqn}) are
$$
   \Delta  =    \frac1{sin\theta}
[\frac{\p }{\p \theta}(sin\theta\frac{\p}{\p \theta})
+\frac1{sin\theta}\frac{\p^2 }{\p \phi^2} ],
$$
$$ 
   A  =  \Delta + \frac{\p }{\p \zeta}(N^2\frac{\p }{\p \zeta}),
$$
$$
  J(p, q)  = \frac{\p p}{\p \theta}\frac1{sin\theta} \frac{\p q}{\p \phi}
	- \frac1{sin\theta}\frac{\p p}{\p \phi} \frac{\p q}{\p \theta},
$$
where $N(\phi, \theta, \zeta) >0 $ is a known smooth function
related to the stratification frequency.
  
The equation (\ref{eqn}) is supplemented by the following
boundary and initial conditions 
\begin{eqnarray}
 \frac{\p\psi}{\p \zeta}  & = & \frac{\p (A\psi)}{\p \zeta} = 0 \quad \mbox{if} \; \zeta =\zeta_0 \; , \label{BC1} \\
 \frac{\p\psi}{\p \zeta} +a \psi & = & \frac{\p (A\psi)}{\p \zeta}+a A\psi =0 \quad \mbox{if} \; \zeta=1 \; , \label{BC2} \\
 \psi |_{t=0}  & = & \psi_0 (\phi, \theta, \zeta)         \;, \label{IC}
\end{eqnarray} 
where $a(\phi, \theta)>0$ is a known smooth function related to heat
transfer between the atmosphere and the earth, and $\psi_0$ is
initial data.

Ddenote $L^2(D)$ as the space of square-integrable functions,  
with the standard norm $\| \cdot \|$.
The problem (\ref{eqn}), (\ref{BC1}, (\ref{BC2}), (\ref{IC})
is well-posed (\cite{Wang}) in  $L^2(D)$, with
solution $\psi(\phi, \theta, \zeta, t)$ at least continuous in time $t$. 

We address the issue of whether there are any  
time-periodic solutions in the nonlinear forced
dissipative quasigeostrophic  dynamics modeled by 
(\ref{eqn}), (\ref{BC1}, (\ref{BC2}). 
To this end we further assume that the forcing 
$f(t): = f(\phi, \theta, \zeta, t)$ is 
periodic in time with period $T>0$, and 
the spatial square-integral of the   forcing $\|f (t)\|$ 
is bounded in time, i.e., $\| f(t) \|$ is bounded by 
a time-independent constant. 
Then we can follow \cite{Wang} exactly to show that 
there is a bounded absorbing set in  $L^2(D)$ (we omit this part).   
That is, all solutions $\psi$ enter a bounded
set $\{ \psi:  \;\;   \|\psi\|   \leq  C(Re, Ro)   \}$
as time goes to infinity. The system (\ref{eqn}), (\ref{BC1}, (\ref{BC2}) 
is  therefore a
dissipative system as defined in \cite{Krasnoselskii}
(also \cite{Temam} or \cite{Hale}).

We now recall a result from \cite{Krasnoselskii}, page 235,
that a $T-$time-periodic nonautonomous dissipative dynamical
system in a Banach space has at least one $T-$time-periodic
solution. This result follows from a Leray-Schauder topological degree
argument and the Browder's principle (\cite{Krasnoselskii}).    
So the system (\ref{eqn}), (\ref{BC1}), (\ref{BC2})
has at least one $T-$time-periodic
solution starting from some initial data $\psi_0$. 
Therefore we   have the following result.

\begin{theorem}
Assume that the forcing $f(\phi, \theta, \zeta, t)$ 
is time-periodic with period $T>0$,
and its spatial square-integral with respect to 
$\phi, \theta, \zeta $ is bounded in time.
 
Then the forced viscous three dimensional baroclinic 
quasigeostrophic model
\begin{eqnarray}
 (A \psi)_t +J(\psi, Ro \; A\psi + 2cos\theta) - \frac1{Re}A^2\psi
	& = & f(\phi, \theta, \zeta, t) ,      \\
 \frac{\p\psi}{\p \zeta} = \frac{\p (A\psi)}{\p \zeta} & = & 0 \quad \mbox{if} \; \zeta =\zeta_0 \; ,   \\
 \frac{\p\psi}{\p \zeta} +a \psi = \frac{\p (A\psi)}{\p \zeta}+a A\psi & = & 0 \quad \mbox{if} \; \zeta=1 \; ,  \\	
\end{eqnarray} 
has at least one time-periodic solution with period $T>0$, for
some square-integrable initial data $\psi_0$.
\end{theorem}

\section{Discussions}

Forced coherent structures in the two dimensional baroclinic model
were studied in \cite{Pierrehumbert}. A wind forced two dimensional baroclinic 
model was also used in the study of multiple geophysical equilibria
(e. g., \cite{Cessi}). In general, it is very difficult to show
existence of periodic coherent structures in spatially
extended physical systems.
In this paper, we have shown that the three dimensional
time-periodic quasigeostrophic patterns may form due to time-periodic 
forcing such as external source of heating.


\begin{thebibliography}{20}

\bibitem{Charney} J. G. Charney, The dynamics of long waves in a baroclinic 
westerly current,
{\em J. Meteorol.} {\bf } (1947), 135-163.


\bibitem{Pedlosky} J. Pedlosky, {\it Geophysical Fluid Dynamics},
Springer-Verlag, 2nd edition, 1987.
 
 
\bibitem{Gill} A. E. Gill, {\em Atmosphere-Ocean Dynamics},
Academic Press, New York, 1982.
 
 
 
\bibitem{Beale} A. J. Bourgeois and J. T. Beale,
Validity of the quasigeostrophic model for large-scale flow in the
atmosphere and ocean, {\em SIAM J. Math. Anal.} {\bf 25} (1994),
1023-1068.
 

\bibitem{Embid_Majda} P. F. Embid and A. J. Majda,
Averaging over fast gravity waves for geophysical flows with arbitrary
potential vorticity, {\em Comm. PDEs} {\bf 21} (1996), 619-658.
  
\bibitem{Desjardins1} B. Desjardins  and E. Grenier,
Derivation of quasigeostrophic potential vorticity equations,
preprint, 1997.

\bibitem{Desjardins2} B. Desjardins  and E. Grenier,
On the homogeneous model of wind driven ocean circulation,
preprint, 1997.

  
  
  
\bibitem{Dutton} J. A. Dutton, The nonlinear quasi-geostrophic equation:
Existence and  uniqueness of solutions on a bounded domain,
{\em J. Atmos. Sci.} {\bf 31} (1974), 422-433.


\bibitem{Bennett} A. F. Bennett and P. E. Kloeden,
The dissipative quasigeostrophic equations,
{\em Mathematika} {\bf 28} (1981), 265-285.
  
  
\bibitem{Holm} D. D. Holm, Hamiltonian formulation of the baroclinic
quasigeostrophic fluid equations, {\em Phys. Fluids} {\bf 29} (1986),
7-8.  

 
 

\bibitem{Wang} S. Wang, Attractors for the 3D baroclinic quasi-geostrophic
equations of large-scale atmosphere, {\em J. Math. Anal. Appl.} {\bf 165}
(1992) 266-283. 
 

\bibitem{Krasnoselskii} M. A. Krasnoselskii and P. P. Zabreiko,
{\em Geometrical Methods of Nonlinear Analysis},
Springer-Verlag, New York, 1984.
 
\bibitem{Temam} R. Temam, {\em  Infinite-Dimensional Dynamical Systems
in Mechanics and Physics,} Springer-Verlag, New York, 1988.
 
 
\bibitem{Hale}  J. K. Hale, {\em Asymptotic Behavior of Dissipative
Systems}, {\rm American Math. Soc.}, 1988.
 
 
\bibitem{Pierrehumbert} R. T. Pierrehumbert and P. Malguzzi,
Forced coherent structures and local multiple equilibria in a
barotropic atmosphere, 
{\em J. Atmos. Sci.} {\bf 41} (1984), 246-257.
 
 
 
\bibitem{Cessi} P. Cessi and G. R. Ierley,
Symmetry-breaking multiple equilibria in quasigeostrophic,
wind-driven flows, {\em J. Phys. Oceanography}, {\bf 25} (1995),
1196-1205. 


\end{thebibliography}
\end{document}